\documentclass[12pt,oneside,english]{amsart}
\usepackage[latin1]{inputenc}
\usepackage{amssymb}

\makeatletter
\newtheorem{theorem}{Theorem}
\newtheorem{lemma}{Lemma}
\newtheorem{example}{Example}

\newtheorem{remark}{Remark}
\newtheorem{conjecture}{Conjecture}

\newtheorem{tab}{Table}

\usepackage{babel}

\makeatother
\begin{document}
\baselineskip=17pt

\title{A conjecture on primes and a step towards justification*}

\author{Vladimir Shevelev}
\address{Departments of Mathematics \\Ben-Gurion University of the
 Negev\\Beer-Sheva 84105, Israel. e-mail:shevelev@bgu.ac.il}

\subjclass{Primary 11N05, Secondary 11B75, 05A19} \keywords{Evil and
odious numbers, primes, counter functions for evil and odious
primes, combinatorial identities. Concerned the sequences: A000069,
A001969, A027697, A027699}

\begin{abstract}
We put a new conjecture on primes from the point of view of its
binary expansions and make a step towards justification.
\end{abstract}

\maketitle
\section{Introduction and main results}\label{s1}
Consider the partition of the set $\mathbb{N}$ into the following
two disjoint subsets
\begin{equation}\label{1}
\mathbb{N}=\mathbb{N}^e\cup \mathbb{N}^o,
\end{equation}
where $\mathbb{N}^e(\mathbb{N}^o)$ is the set of positive integers
which have even (odd) number of 1`s in their binary expansions.
These numbers are called the evil and the odious numbers
respectively \cite{Slo-1}. There are some results for these numbers
and some applications of them in
\cite{Allouche-1992},\cite{Allouche-2003},\cite{Guy-1991},\cite{Guy-1995},
\cite{LamMo},\cite{Mc1}.

Consider the same partition of the set $\mathbb{P}$ of prime numbers
\cite{Slo-2}:

\begin{equation}\label{6}
\mathbb{P}=\mathbb{P}^e\cup \mathbb{P}^o.
\end{equation}
For example, all the Fermat primes are evil while all the Mersenne
primes $> 3$ are odious.

     Using direct calculations up to $10^9$ we noticed that among the
primes not exceeding $n$ the evil primes are never in majority
except for the cases $n=5 \; and \; n=6$.  Moreover, in the
considered limits the excess of the odious primes is not monotone
but increases on the whole with records on primes
$2,13,41,67,79,109,131,137,\ldots$

Let $\pi^e(x)(\pi^o(x))$ denote the number of the evil (odious)
primes not exceeding $x$. Put
$$
m_n=\min\limits_{x\in(2^{n-1},2^n)}(\pi^o(x)-\pi^e(x)).
$$

The following table shows that $m_n$ increases monotonically.

\begin{tab}\label{tab1}

\begin{small}
$$
\begin{matrix}
&n &m_n \;\;\;&n &m_n\\[8pt]
&5 &0 \;\;\;&19  &1353\\
&6 &2\;\; \;&20  &1855\\
&7 &4\;\; \;&21  &3659\\
&8 &7\;\; \;&22  &5221\\
&9 &13\;\; \;&23  &10484\\
&10 &19\;\; \;&24  &14933\\
&11 &39\;\; \;&25  &27491\\
&12 &54\;\; \;&26  &35474\\
&13 &104\;\; \;&27  &68816\\
&14 &139\;\; \;&28  &97342\\
&15 &251\;\; \;&29  &186405\\
&16 &334\;\;\; &30 &265255\\
&17 &590\\
&18 &716
\end{matrix}
$$
\end{small}
\end{tab}

Therefore, the following conjecture seems plausible.

\begin{conjecture}\label{conj1}   For all
$n\in\mathbb{N},\;n\neq 5,6 $\upshape
\end{conjecture}
\begin{equation}
\pi ^e(n)\leq\pi ^o(n);
\end{equation}

\itshape moreover,\upshape
\begin{equation}\label{100}
\lim_{n\rightarrow\infty}(\pi^o(n)-\pi^e(n))=+\infty.
\end{equation}

For $a$ positive integer $a$, denote $\mu^e_a(n)(\mu_a^o(n))$ the
number of \itshape odd \upshape evil (odious) nonnegative integers
divisible by $a$ and less than $n$.

\begin{remark}   We include in this definition 0
(which is an evil integer) and use "less than" instead of "not
exceeding" for the sake of more simplicity of the formulas which
appear below.
\end{remark}
Put
\begin{equation}\label{101}
\Delta^{odd}_a(n)=\mu^e_a(n)-\mu^o_a(n)
\end{equation}

\begin{theorem}\label{thm1}\;\;Let $p,q,\ldots$ denote odd primes. Then

\begin{equation}\label{11}
\pi^o(n)-\pi^e(n)=\varepsilon_n+\sum_{p\leq
n}\Delta_p^{odd}(n)-\sum_{p< q\leq n}\Delta_{p,q}^{odd}(n)+\ldots,
\end{equation}

where  $|\varepsilon_n|\leq 4$.
\end{theorem}
In this article we make only the first step of investigation of
$\pi^o(n)-\pi^e(n)$ with help of (\ref{11}). Namely, by
combinatorial methods we study in detail $\Delta_3^{odd}(n)$.

Let

\begin{equation}\label{16}
\Delta_3^{odd}([a,b))=\Delta_3^{odd}(b)-\Delta_3^{odd}(a)
\end{equation}

\begin{theorem}\label{thm2}\;\;1)  $\Delta^{(odd)}_3([0,2^n))=
3^{\lfloor\frac n 2\rfloor -1},\;n\geq 2$

2)$\Delta^{(odd)}_3([2^n,2^n+2^m))=\begin{cases}0,\;
n\;and\;m\;are\;even,\;\; 2\leq m\leq n-2,\\
3^{\frac{m-2}{2}},\;n \;is \;odd, \;m \;is \;even,\;\;2\leq m\leq n-1,\\
-3^{\frac{m-3}{2}},\;n\;is\;even,\;m\; is\; odd,\;\; 3\leq m\leq n-1,\\
2\cdot 3^{\frac{m-3}{2}},n\;and\;m\;are\;odd,\;3\leq m\leq n-2.
\end{cases}$

\; \; \;\;
 3)$\Delta^{(odd)}_3([2^n +2^{n-2},
2^n+2^{n-2}+2^m))=\begin{cases}-3^{\frac{m-2}{2}},\;
n\;and\;m\;are\;even,\;\;2\leq m\leq n-4,\\
0,\;n \;is \;odd, \;m \;is \;even,\;\;2 \leq m\leq n-3, \\
-2\cdot 3^{\frac{m-3}{2}},\;n\;is\;even,\;m\; is\; odd, \;\; 3\leq m\leq n-3,\\
3^{\frac{m-3}{2}},n\;and\;m\;are\;odd,\;3\leq m\leq n-4.
\end{cases}$
\end{theorem}

Consider together with $\Delta^{(odd)}_3([a,b))$ also
$\Delta^{(even)}_3([a,b))$ which means the difference between the
numbers of evil and odious \itshape even \upshape integers divisible
by 3 on $[a,b)$. Put

\begin{equation}\label{17}
\Delta_3([a,b))=\Delta^{(odd)}_3([a,b))+\Delta^{(even)}_3([a,b))
\end{equation}

\begin{theorem}\label{thm3}

1)$\Delta_3([0, 2^n ))=\begin{cases}2\cdot 3^{\frac n 2 -1},\;
n\;is \;even\\
3^{\frac{n-1}{2}},\;n\;is\; odd,\;n\geq 1
\end{cases}$

\; \; \; \;

2)$\Delta_3([2^n,2^n+2^m))=\begin{cases}3^{\lfloor
\frac{m-1}{2}\rfloor},\;if \; n \; is \; odd, 1\leq m \leq {n-1},\\
3^{\frac m 2-1},\; if \; n \;and \;m \;are \;even,\;2\leq m\leq n-2, \\
0,\; if \;n\;is\;even,\;m\; is\; odd, 1\leq m \leq n-1,
\end{cases}$

\; \; \; \;

3)$\Delta_3([2^n+2^{n-2}, 2^n+2^{n-2}+2^m))=\begin{cases}-3^{\lfloor
\frac{m-1}{2}\rfloor},\;if \; n \; is \; even,1\leq m \leq n-3,\\
-3^{\frac m 2-1},\; if \; n \;is \;odd, \;m \; is \;even,\;2\leq  m\leq n-3, \\
0,\; if \;n\;and\;m\; are\; odd, 1\leq m\leq n-4
\end{cases}$
\end{theorem}
 At last, the following result is valid.

\begin{theorem}\label{thm4}
$$
\lim_{n\rightarrow\infty}\frac{\ln \Delta_3^{(odd)}([0,n))}
{\ln n}=\frac{\ln 3}{\ln 4}.
$$
\end{theorem}

Using theorem \ref{thm4} and simple heuristic arguments we put our
conjecture in the following quantitative form.

\begin{conjecture}
$$
\lim_{n\rightarrow\infty}\frac{\ln(\pi^o(n)-\pi^e(n))}{\ln
n}=\frac{\ln 3}{\ln 4}.
$$
\end{conjecture}

Conjecture 2 is illustrated by Table 2 in Section 3.

     In the following Section we prove Theorems 1-4. Section 3 is
devoted to some heuristic arguments which lead to Conjecture 2.
Finally, in Section 4 we consider the increment of the excess of
odiores  primes on$(0,2^n) (Table 3)$.

\section{Proofs of results}\label{s2}

\bfseries   A. Proof of Theorem 1.\mdseries

Denote $\nu^e(n)(\nu^o(n))$ the number of evil (odious) nonnegative
integers on interval $[0,n)$.

\begin{lemma}\label{lem1}

We have
\begin{equation}\label{102}
|\nu^o(n)-\nu^e(n)|\leq 1, \;\;n\in \mathbb{N}.
\end{equation}
\end{lemma}

Proof. The Lemma follows from the identity
\begin{equation}\label{103}
\nu^e(2m)=\nu^o(2m), \;\;m\in \mathbb{N},
\end{equation}
which is proved by induction.

     Notice that (\ref{103}) is satisfied for $m=1$. Assuming that
it is valid for $2m$  we prove (\ref{103}) for $2(m+1)$. Indeed, let
$m$ has $k$ 1`s in the binary expansion. Then we have evidently

$$
\nu^e(2m+1)-\nu^o(2m+1)=(-1)^k.
$$

On the other hand, the last number in interval $[0,2m+2)$, i.e. the
number $2m+1$ has  $k+1$ 1`s and thus
$\nu^e(2m+2)-\nu^o(2m+2)=0.\;\;\blacksquare$

     Let $\lambda^e(n)(\lambda^o(n))$ denote the number of \itshape
even \upshape evil (odious) numbers less than $n$. At last, denote
$\sigma^e(n)(\sigma^o(n))$ the number of evil (odious) \itshape odd
composite \upshape numbers less than $n$.

 For $n\geq 3$ we have
\begin{equation}\label{104}
\pi^o(n)-\pi^e(n)+\sigma^o(n)-\sigma^e(n)+\lambda^o(n)-\lambda^e(n)-1=\nu^o(n)-
\nu^e(n)+\delta_n,
\end{equation}

where according to the definition of $\pi^o(n)(\pi^e(n)), \delta_n$
is 1, if $n$ is an odious prime, -1, if $n$ is an evil prime,
0-otherwise. Subtraction 1 in the left hand side of (\ref{104})  is
connected with the fact that only 2 is an odious prime and
simultaneously is an odious even integer.

Using Lemma 1 and the evident identity

\begin{equation}\label{105}
\lambda^o(n)-\lambda^e(n)=\nu^o(\frac n 2)-\nu^e(\frac n 2)
\end{equation}

we find from (\ref{104})
\begin{equation}\label{106}
\pi^o(n)-\pi^e(n)=\sigma^e(n)-\sigma^o(n)+\varepsilon_m,
\end{equation}

where $|\varepsilon_n|\leq 4$.

At last, by inclusion-exclusion from (\ref{106}) we obtain
(\ref{11}) $\blacksquare$

\bfseries  B.Proofs of Theorems \ref{thm2}-\ref{thm3}.\mdseries \;\;
It is easy to see that for nonnegative integers $a<b$

\begin{equation}\label{18}
\Delta^{(even)}_3([2a, 2b))=\Delta_3([a, b))
\end{equation}

and consequently
\begin{equation}\label{19}
\Delta^{(odd)}_3([2a, 2b))=\Delta_3([2a, 2b))-\Delta_3([a, b)).
\end{equation}

Therefore, it is sufficient to prove Theorem \ref{thm3}  and by
(\ref{18})-(\ref{19}) we shall get also Theorem \ref{thm2}. For the
proof of Theorem \ref{thm3} we need a simple lemma.

\begin{lemma}\label{lem2}

Let for a nonnegative integer $n,\;i^{even}(n)(i^{odd}(n))$ denote
the number of even (odd) powers of 2 in the binary representation of
$n$. Then
\begin{equation}\label{20}
n\equiv 0(mod 3)\Leftrightarrow i^{even}(n)\equiv i^{odd}(n)(mod 3)
\end{equation}
\end{lemma}
\bfseries     Proof. 1.  Straightforward.\mdseries

\bfseries Proof of Theorem \ref{thm3}.\mdseries

1a) let $n$ be even, $n=2m$. Consider all the nonnegative integers
not exceeding \;\;$2 ^{2m}-1$ \;\;which have\;\; $2m $\;\; binary
positions with numbering \;\;\;\;\;\;$0,1,\ldots,2m-1$ beginning
from the right. To find the difference between the numbers of evil
and odious integers divisible by 3 not exceeding $2^{2m} -1$, let
choose j even position for 1`s (and $m-j$ even position for 0`s) and
according to Lemma \ref{lem2} let choose $j+3k\;(k\geq 0)$ odd
position for 1`s (and the rest of the odd positions for 0`s).

After that, vice versa, we choose $j$ odd positions for 1`s (and
$n-j$ odd positions for 0`s) and $j+3k\;(k\geq 1)$ even positions
for 1`s (and the rest of the even positions for 0`s). Notice that,
for each $j$ the parity of the number of the chosen  1`s is the same
as the parity of k. Thus

\begin{equation}\label{21}
\Delta_3([0, 2^{2m}))=\sum_{j\geq 0}\binom{m}{j}^2+2\sum_{k\geq
1}(-1)^k \sum_{j\geq 0}\binom{m}{j}\binom{m}{j+3k}.
\end{equation}

Since (cf.\cite{Rio},p.8)

\begin{equation}\label{22}
\sum_{j\geq 0}\binom{m}{j}\binom{m}{j+3k}=\binom{2m}{m+3k},\;k\geq
0,
\end{equation}

then by (\ref{21})
\begin{equation}\label{23}
\Delta_3([0,2^{2m}))=\binom{2m}{m}+2\sum_{k\geq 1}(-1)^k
\binom{2m}{m+3k}.
\end{equation}

To calculate $\sum_{k\geq 1}(-1)^k \binom{2m}{m+3k}$ in (\ref{23})
we need some lemmas.

\begin{lemma}\label{lem3}\;\;(\cite{Rio},p.86)

\begin{equation}\label{24}
\sum_{k=0}^{\lfloor\frac m 3\rfloor}\binom{2m}{m+3k}=\frac 1 3
(2^{2m-1}+1)+\frac 1 2 \binom{2m}{m}.
\end{equation}
\end{lemma}

\begin{lemma}\label{lem4}
\begin{equation}\label{25}
\sum_{t=0}^{\lfloor\frac m 6\rfloor}\binom{2m}{m+6t}=\frac 1 2
\left(\frac{2^{2m-1}+1}{3}+3^{m-1}+ \binom{2m}{m}\right).
\end{equation}
\end{lemma}

Proof. Denote the left hand side of (\ref {25}) by $\sigma(m)$.

Let  $m=6l+s, \;\;0\leq s\leq 5$. Then

\begin{equation}\label{26}
\sigma(m)=\sum^l_{t=0}\binom{12l+2s}{6l+s-6t}=\sum^l_{k=0}\binom{12l+2s}{6k+s}
\end{equation}

Together with $\sigma(m)$ let consider the sum
$$
\sigma_1(m)=\Sigma^{2l}_{k=l+1}\binom{12l+2s}{6k+s}=(2l-k=t)=
$$
\begin{equation}\label{27}
=\Sigma^{l-1}_{t=0}
\binom{12l+2s}{12l-6t+s}=\Sigma^{l-1}_{t=0}\binom{12l+2s}{6t+s},
\end{equation}

From  (\ref {26}),(\ref {27}) we conclude that

\begin{equation}\label{28}
\sigma(m)=\sigma_1(m)+\binom{12l+2s}{6l+s},
\end{equation}

Consequently,
$$
2\sigma(m)=\sigma(m)+\sigma_1(m)+\binom{12l+2s}{6l+s}=
\sum^{2l}_{k=0}\binom{2m}{6k+s}+\binom{2m}{m}.
$$
Thus,
\begin{equation}\label{29}
\sum^{\lfloor\frac m 6\rfloor}_{t=0}\binom{2m}{m+6t}=\frac 1 2
\left(\sum_{k=0}^{\frac{m-s}{3}}\binom{2m}{6k+s}+\binom{2m}{m}\right),
\end{equation}

where $0\leq s \leq 5$.

Notice that, $\frac{m-s}{3}$ is the "natural" upper limit of the sum
on the right hand side in (\ref{29}).Indeed, in this sum $k\leq
\lfloor\frac{2m-s}{6}\rfloor=\lfloor\frac{12l+s}{6}\rfloor=2l=\frac{m-s}{3}$.
 To calculate this sum we use the formula (\cite {Rio}, p.161)from which for $s=m-6t$ it
follows that
$$
\sum^{\frac{m-s}{3}}_{k=0}\binom{2m}{6k+s}=\frac 1
6\sum^6_{j=1}e^{\frac{\pi i}{3}(-jm)}\left( 1+e^{\frac{\pi
i}{3}j}\right)^{2m}=
$$
$$
=\frac 1 6\left(e^{-\frac{\pi m}{3}i}\left(1+e^{\frac{\pi
i}{3}}\right)^{2m}+e^{-\frac{2\pi m}{3}i}\left(1+e^{\frac{2\pi
i}{3}}\right)^{2m}+e^{-\frac{4\pi m}{3}i}\left(1+e^{\frac{4\pi
i}{3}}\right)^{2m}+\right.
$$
$$
+e^{-\frac{5\pi m}{3}i}\left(1+e^{\frac{5\pi
i}{3}}\right)^{2m}+2^{2m}\left.\right)= \frac 1 6
\left(\right.e^{-\frac{\pi m}{3}i}\left(1+e^{\frac{\pi
i}{3}}\right)^{2m}+ e^{-\frac{2\pi m}{3}i}\left(1+e^{\frac{2\pi
i}{3}}\right)^{2m}+
$$

$$
+e^{\frac{2\pi m}{3}i}\left(1+e^{-\frac{2\pi
i}{3}}\right)^{2m}+e^{\frac{\pi m}{3}i}\left(1+e^{-\frac{\pi
i}{3}}\right)^{2m}+2^{2m}\left.\right)=\frac 1 3
\left(2^{2m-1}+\right.
$$
$$
+Re\left(e^{-\frac{\pi m}{3}i}\left(1+e^{\frac{\pi
i}{3}}\right)^{2m}\right)+Re\left(e^{-\frac{2\pi
m}{3}i}\left(1+e^{\frac{2\pi i}{3}}\right)^{2m}\right)\left.\right).
$$

Noticing that, $1+e^{\frac{\pi i}{3}}=\frac 3 2
+\frac{\sqrt{3}}{2}i=\sqrt {3} e^{\frac{\pi i}{6}},\;1+e^{\frac{2\pi
i}{3}}=e^{\frac \pi 3 i}$ we have

\begin{equation}\label{30}
\Sigma^{\frac{m-s}{3}}_{k=0}\binom{2m}{6k+s}=\frac 1 3
\left(2^{2m-1}+3^m+1\right)
\end{equation}

and by (\ref{29}), (\ref{30}) we obtain the lemma $\blacksquare$

\begin{lemma}\label{lem5}
\begin{equation}\label{31}
\sum_{k\geq 0}(-1)^k\binom{2m}{m+3k}=3^{m-1}+\frac 1 2
\binom{2m}{m}.
\end{equation}
\end{lemma}

Proof.  We have
$$
\sum_{k\geq 0}(-1)^k\binom{2m}{m-3k}+\sum_{k\geq
0}\binom{2m}{m-3k}=2\sum_{j\geq 0}\binom{2m}{m-6j}
$$
and by Lemmas \ref{lem3}, \ref{lem4} we obtain the lemma
$\blacksquare$

Now from  (\ref{23}) and Lemma \ref{lem5} we find

$$
\Delta_3\left([0,2^{2m})\right)=2 \cdot 2^{m-1}.
$$

1b) As opposed to the case 1a) here we have $2m-1$ positions from
which $m$ even and $m-1$ odd. Hence, by the same combinatorial
arguments we find
$$
\Delta_3\left([0,2^{2m-1})\right)=\sum_{j\geq
0}\binom{m}{j}\binom{m-1}{j}+
$$
\begin{equation}\label{32}
+\sum_{k\geq 1}(-1)^k\left(\sum_{j\geq
0}\binom{m}{j}\binom{m-1}{j+3k}+\sum_{j\geq
0}\binom{m-1}{j}\binom{m}{j+3k}\right)
\end{equation}

Since (cf.\cite{Rio},p.8)

$$
\sum_{j\geq 0}\binom{m}{j}\binom{m-1}{j+3k}=\binom{2m-1}{m+3k}
$$
$$
\sum_{j\geq 0}\binom{m-1}{j}\binom{m}{j+3k}=\binom{2m-1}{m+3k-1}
$$

then by (\ref{32}) and Lemma \ref{lem5} we have

$$
\Delta_3\left([0,2^{2m-1})\right)=\binom{2m-1}{m}+\sum_{k\geq
1}(-1)^k\binom{2m}{m+3k}=
$$
$$
=-\frac 1 2 \binom{2m}{m}+\sum_{k\geq
0}(-1)^k\binom{2m}{m+3k}=3^{m-1}\;\;\blacksquare
$$

2a) Let $m$ be even, $m=2l$.  Let, for definiteness, $n$ be even.
Choose $j-1$ of the last $l$ even positions for 1`s (and the rest
$l-(j-1)$ positions for 0`s) and according to Lemma \ref{lem2},
choose $j+3k \;(k\geq 0)$ of the last $l$ odd positions for 1`s (and
$l-j-3k$ positions for 0`s). After that, vice versa, we choose $j$
of the last $l$ odd positions for 1`s and also $j-1+3k$ of the last
$l$ even positions for 1`s and the rest of the positions for 0`s.
For each $j$ the parity of the number of all 1`s (including the 1
corresponding to $2^n$) is the same as the parity of $k$. Thus,

$$
\Delta_3\left([2^n,2^n+2^{2l})\right)=\sum_{j\geq 1}\binom
{l}{j-1}\binom{l}{j}+
$$
\begin{equation}\label{33}
+\sum_{k\geq 1}(-1)\left(\sum_{j\geq
1}\binom{l}{j-1}\binom{l}{j+3k}+ \sum_{j\geq
0}\binom{l}{j}\binom{l}{j-1+3k} \right).
\end{equation}

Since
$$
\sum_{j\geq 1}\binom{l}{j-1}\binom{l}{j+3k}=\binom{2l}{l+3k+1},
$$
$$
\sum_{j\geq 0}\binom{l}{j}\binom{l}{j+3k-1}=\binom{2l}{l+3k-1},
$$

then
$$
\Delta_3\left([2^n,2^n+2^{2l})\right)=\binom{2l}{l+1}+
$$
\begin{equation}\label{34}
+\sum_{k\geq
1}(-1)^k\left(\binom{2l}{l+3k-1}+\binom{2l}{l+3k+1}\right)
\end{equation}

It is easy to verity that
\begin{equation}\label{35}
\binom{2l}{l+3k-1}+\binom{2l}{l+3k+1}=\binom{2l+2}{l+3k+1}-2\binom{2l}{l+3k}.
\end{equation}

Thus, using Lemma \ref{lem5} for $m=l$ and $m=l+1$ we have
$$
\Delta_3[2^n,2^n+2^{2l})=\binom{2l}{l+1}+3^l-\frac 1 2
\binom{2(l+1)}{l+1}-2\cdot 3^{l-1}+\binom{2l}{l}=
$$

\begin{equation}\label{36}
=\binom{2l}{l+1}-\binom{2l+1}{l+1}+\binom{2l}{l}+3^{l-1}=3^{l-1}.
\end{equation}

It is evident that in this case the validity of (\ref{33}) does not
depend on the parity of $n$.

2b) Let $m$ be odd, $m=2l+1, \;\;l\geq 0$. As opposed to the case 2a
here we have the last $2l+1$ positions from which $l+1$ are even and
$l$ are odd. Hence, by the same arguments we find
$$
\Delta_3\left([2^n, 2^n+2^{2l-1}) \right)=\sum_{j\geq
1}\binom{l}{j}\binom{l+1}{j-1}+\sum_{k\geq 1}(-1)^k\left(
\sum_{j\geq 0}\binom{l}{j+3k+1}\binom{l+1}{j}\right.+
$$
\begin{equation}\label{37}
+\left.\sum_{j\geq 0}\binom{l+1}{j+3k-1}\binom{l}{j}\right),
\;\;if\; n\;is\;even,
\end{equation}

and

$$
\Delta_3\left([2^n, 2^n+2^{2l-1}) \right)=\sum_{j\geq
1}\binom{l}{j-1}\binom{l+1}{j}+\sum_{j\geq 1}(-1)^k\left(
\sum_{j\geq 1}\binom{l+1}{j+3k}\binom{l}{j-1}\right.+
$$

\begin{equation}\label{38}
+\left.\sum_{j\geq 0}\binom{l}{j+3k-1}\binom{l+1}{j}\right),
\;\;if\; n\;is\;odd.
\end{equation}

Now by (\ref{37}) for even $n$ we have
$$
\Delta_3\left([2^n, 2^n+2^{2l-1}) \right)=\sum_{j\geq
0}\binom{l+1}{j}\binom{l}{j+1}+
$$
$$
+\sum_{k\geq 1}(-1)^k\left(
\binom{2l+1}{l+3k+2}+\binom{2l+1}{l+3k-1}\right)=
$$
$$
=\binom{2l+1}{l+2}+\sum_{k\geq 1}(-1)^k
\binom{2l+1}{l+3k+2}+\sum_{k\geq 1}(-1)^k
\binom{2l+1}{l+3k-1}=
$$
$$
=\binom{2l+1}{l+2}+\sum_{k\geq 1}(-1)
\binom{2l+1}{l+3k+2}-\sum_{k\geq 0}(-1)^k \binom{2l+1}{l+3k+2}=
$$
$$
=\binom{2l+1}{l+2}-\binom{2l+1}{l+2}=0.
$$

and by (\ref{38}) for odd $n$ we have

$$
\Delta_3\left([2^n, 2^n+2^{2l-1}) \right)=\sum_{j\geq
0}\binom{l}{j}\binom{l+1}{j+1}+\sum_{k\geq 1}(-1)^k
\left(\sum_{j\geq 0}\binom{l}{j}\binom{l+1}{j+3k+1}+\right.
$$
$$
+\left.\sum_{j\geq 0}\binom{l+1}{j}\binom{l}{j+3k-1}\right)=
\binom{2l+1}{l+1}+\sum_{k\geq
1}(-1)^k\left(\binom{2l+1}{l+3k+1}+\right.
$$
$$
+\left.\binom{2l+1}{l+3k}\right)=\binom{2l+1}{l+1}+\sum_{k\geq
1}(-1)^k\binom{2l+2}{l+3k+1}=
$$
$$
=\binom{2l+1}{l+1}+\sum_{k\geq 1}(-1)^k\binom{2(l+1)}{(l+1)-3k},
$$

and by Lemma \ref{lem5} for odd $n$ we obtain
$$
\Delta_3\left([2^n, 2^n+2^{2l-1})
\right)=\binom{2l+1}{l+1}+3^l-\frac 1 2 \binom{2l+2}{l+1}=3^l\;\;
\blacksquare
$$

3) denote by $\Delta_{3,h}([a,b)),\; (h\in \mathbb{N})$, the
difference between the numbers of evil and odious integers on
$[a,b)$ having the form $3t+i,\;\;i=1,2$, where $i\equiv h(\mod{3})$

\begin{lemma}\label{lem6}

1) $\Delta_{3,1}([0,2^n))=\begin{cases}-3^{\frac n 2 -1},\; if\; n
\;is \;even\\0,\; if \; n\; is\; odd \end {cases}$

2)$\Delta_{3,2}([0,2^n))= -3^{\lfloor\frac{n-1}{2}\rfloor}$

\end{lemma}

Proof. Notice that,
\begin{equation}\label{39}
\Delta_3([2^n,2^n+2^m))=-\Delta_{3,2^n}([0,2^m))
\end{equation}

Since by $\mod{3}$
$$
2^n\equiv \begin{cases} 1,\; if \; n \; is \; even\\2,\; if \; n \;
is \; odd,\end{cases}
$$
then by (\ref{39})

$$
\Delta_{3,1}([0,2^m))=\Delta_3([2^n, 2^n+2^m)) if \; n \; is\; even
$$
$$
\Delta_{3,2}([0,2^m))=\Delta_3([2^n, 2^n+2^m)) if \; n \; is\; odd
$$

and the lemma follows from the previous point  $\blacksquare$.

Now we are able to complete the proof of Theorem \ref{thm3}.

a) Let $n$ be even, $n=2t$.

We have

\begin{equation}\label{40}
\Delta_3
([2^{2t}+2^{2t-2},\;2^{2t}+2^{2t-2}+2m))=\Delta_{3,2^{2t}+2^{2t-2}}([0,2^m)).
\end{equation}

Since
$$
2^{2t}+2^{2t-2}\equiv 5\cdot 2^{2t-2}\equiv 2 \;\;(\mod{3}),
$$

then by (\ref{40}) and by Lemma \ref{lem6}

$$
\Delta_3
([2^{2t}+2^{2t-2},\;2^{2t}+2^{2t-2}+2^m))=\Delta_{3,2}([0,2^m))=-3^{\lfloor
\frac {m-1}{2}\rfloor}.
$$

b) Let now $n$ be odd, $n=2t+1$. Since
$$
2^{2t+1}+2^{2t-1}\equiv 5 \cdot 2^{2t-1}\equiv 1\;\;(\mod{3})
$$

then using Lemma \ref{lem6} we have
$$
\Delta_3
([2^{2t+1}+2^{2t-1},\;2^{2t+1}+2^{2t-1}+2^m))=\Delta_{3,1}([0,2^m))=\begin{cases}
-3^{\frac m 2 -1}, \;if\; m\; is \;even\\0,\;if\; m \; is \;
odd\end{cases}.
$$

This completes the proof of both Theorem \ref{thm3} and, in view of
(\ref{19}), Theorem \ref{thm2}\; $\blacksquare$

Notice that,the results of Theorems \ref{thm2},\ref{thm3} one can
write in terms of the counting functions of the corresponding
sequences. For example, let us consider the first points of these
theorems. Let $\nu^e_3(n)(\nu^o_3(n))$  denote the number of the
evil (odious) divisible by 3 nonnegative integers less than $n$.
Then from the first point of Theorem \ref{thm3} for $n\geq 1$ we
have
$$
\nu^e_3(2^n)=\begin{cases}\frac 1 2
\left(\frac{2^n+1}{3}+3^{\frac{n-1}{2}}\right),\;\;if
\;\;n\;\;is\;\;odd\\ \frac{2^{n-1}+1}{3}+3^{\frac n 2 -1},\;\;if
\;\;n\;\;is\;\;even;
\end{cases}
$$

$$
\nu^o_3(2^n)=\begin{cases}\frac 1 2
\left(\frac{2^n+1}{3}-3^{\frac{n-1}{2}}\right),\;\;if
\;\;n\;\;is\;\;odd\\ \frac{2^{n-1}+1}{3}-3^{\frac n 2 -1},\;\;if
\;\;n\;\;is\;\;even.
\end{cases}
$$

     Furthermore, let as above $\mu^e_3(n)(\mu^{(o)}_3(n))$ denote
the number of the evil (odious) divisible by 3 nonnegative \itshape
odd \upshape integers less than $n$. Then from the first point of
Theorem \ref{thm2} for $n\geq 2$ we have

$$
\mu^e_3(2^n)=\frac 1 2
\left(\left\lfloor\frac{2^{n-1}+1}{3}\right\rfloor+3^{\left\lfloor\frac
n 2\right\rfloor-1}\right),
$$

$$
\mu^o_3(2^n)=\frac 1 2
\left(\left\lfloor\frac{2^{n-1}+1}{3}\right\rfloor-3^{\left\lfloor\frac
n 2\right\rfloor-1}\right).
$$

Notice in addition that, Theorem \ref{thm2} (Theorem \ref{thm3})
allows to calculate for any $n$ the number
$\Delta_3^{(odd)}([0,n))(\Delta_3([0,n)))$.

Indeed, let
\begin{equation}\label{41}
\Delta_3^{(odd)}=\Delta_3^{(odd)}(2^{n_1}+2^{n_2}+ \ldots
+2^{n_k},\; 2^{n_1}+2^{n_2}+\ldots +2^{n_k}+2^m).
\end{equation}

Consider the sums
$$
a=\sum_{i:n_i\equiv 0(\mod{2})}1,\;\; b=\sum_{i:n_i\equiv
1(\mod{2})}1.
$$
Let
$$
a\equiv \alpha(\mod{3}),\; b\equiv \beta(\mod{3}),
$$
so that $0\leq \alpha,\beta\leq 2$. Then for any integer $t>\frac m
2$ we have

\begin{equation}\label{42}
\Delta_3^{odd}=\begin{cases}(-1)^k\Delta_3^{odd}([0, 2^m)), \;
if\;\;\alpha=\beta\\(-1)^{k-1}\Delta_3^{odd}([2^{2t},2^{2t}+ 2^m)),
\;if\;\;\alpha-\beta=1,\\(-1)^{k-1}\Delta_3^{odd}([2^{2t+1},2^{2t+1}+
2^m)),\;if\;\;\alpha-\beta=-1,\\(-1)^k\Delta_3^{odd}([2^{2t+2}+2^{2t},2^{2t+2}+2^{2t}+
2^m)),
\;if\;\;\alpha-\beta=2,\\(-1)^k\Delta_3^{odd}([2^{2t+3}+2^{2t+1},2^{2t+3}
+2^{2t+1}+2^m)), \;if\;\;\alpha-\beta=-2.\end{cases}
\end{equation}

(\ref{42}) follows immediately from Lemma \ref{lem2}. The analogous
equality is valid for $\Delta_3$.

\begin{example}   $n=105$. The interval [0, 105)
contains 17 odd numbers divisible by 3. Among them there are 5
odious numbers (namely, 21, 69, 81, 87, 93) and 12 evil numbers.
Thus, $\Delta_3^{odd}([0,105))=7$.
\end{example}

Let find now this value by the algorithm. We have

$$
[0,105)=[0,2^6)\cup [2^6,2^6+2^5)\cup[2^6+2^5,
$$
\begin{equation}\label{43}
2^6+2^5+2^3)\cup[2^6+2^5+2^3,2^6+2^5+2^3+1).
\end{equation}

The last subset does not contain any odd number.  By (\ref{42}) we
have
$$
\Delta_3^{odd}([2^6+2^5,2^6+2^5+2^3))=
$$
\begin{equation}\label{44}
=\Delta_3^{odd}([0,2^3)) \;\;(here \; k=2, \;\alpha=\beta=1).
\end{equation}

Therefore, by (\ref{43}),(\ref{44}) and Theorem \ref{thm2} we find
$$
\Delta_3^{odd}([0,105))=3^2-3+1=7  \; \blacksquare
$$

\bfseries  C.Proof of Theorem \ref{thm4}.\mdseries

In view of (\ref{41})-(\ref{42}) it is sufficient to prove Theorem
\ref{thm4} for the numbers of the form
$$
a) 2^n,\;b) 2^n+2^m, \;\; m\leq n-1,\;c)2^n+2^{n-2}+2^m,\;\;m\leq
n-3.
$$

a) According to the point 1 of Theorem 2 we have

$$
\lim_{n\rightarrow\infty}\frac{\ln \Delta_3^{odd}([0,2^n))}{\ln
2^n}=\lim_{n\rightarrow\infty}\frac{(\lfloor\frac n 2\rfloor-1)\ln
3}{n\ln 2}=\frac{\ln 3}{\ln 4}.
$$

b) According to the points 1 and 2 of Theorem \ref{thm2} and taking
into account that $m\in[2, n-1]$ we have

$$
\Delta_3^{odd}([0, 2^n+2^m))=\Delta_3^{odd}([0,
2^n))+\Delta_3^{odd}([2^n, 2^n+2^m))\leq
$$
$$
\leq 3^{\lfloor\frac n 2\rfloor-1}+2\cdot 3^{\frac{m-3}{2}}\leq
3^{\frac{n-2}{2}}+2\cdot 3^{\frac{n-4}{2}}\leq c_1\cdot 3^{\frac n
2}.
$$

Therefore,
$$
\limsup_{n\rightarrow \infty}\frac{\ln\Delta_3^{odd}([0,
2^n+2^m))}{\ln(2^n+2^m)}\leq \lim_{n\rightarrow \infty}\frac{\ln
c_1+\frac n 2 \ln 3}{n\ln 2}=\frac{\ln 3}{\ln 4}.
$$

On the other hand,
$$
\Delta_3^{odd}([0, 2^n+2^m))\geq 3^{\frac{n-1}{2} -1}-
3^{\frac{m-3}{2}}\geq 3^{\frac{n-3}{2}} - 3^{\frac n 2-2}\geq
0.08\cdot 3^{\frac n 2}.
$$

Thus,
$$
\liminf_{n\rightarrow \infty}\frac{\ln\Delta_3^{odd}([0,
2^n+2^m))}{\ln(2^n+2^m)}\geq \lim_{n\rightarrow \infty}\frac{\ln
0.08 +\frac n 2 \ln 3}{\ln 2+n\ln 2}=\frac{\ln 3}{\ln 4}.
$$

c) Analogously, according to the points 1, 2 and 3 of Theorem
\ref{thm2} and taking into account that $m\in [2, n-3]$ we have

$$
\Delta_3^{odd}([0, 2^n+2^{n-2}+2^m))\leq (c_1+1)\cdot 3^{\frac n 2}+
3^{\frac{n-6}{2}}=c_2 3^{\frac n 2},
$$
$$
\Delta_3^{odd}([0, 2^n+2^{n-2}+2^m))\geq 0.08\cdot 3^{\frac n
2}-2\cdot 3^{\frac{n-6}{2}}\geq 0.005\cdot 3^{\frac n 2}
$$
and we are done . $\blacksquare$

\section{On Conjecture 2}

Show that Conjecture 2 is a corollary of the following heuristic
argument: the behavior of primes with the point of view the excess
of the odious primes is proportionally similar to behavior of
numbers not divisible by 2 and 3. Indeed, the number of the latter
numbers less than $n$ is $n-1-\lfloor\frac {n-1}
{2}\rfloor-\lfloor\frac{n-1}{3}\rfloor +\lfloor\frac
{n-1}{6}\rfloor\sim \frac n 3$. Thus, the excess $\delta(n)$ of the
odious numbers not divisible by 2 and 3 and less than $n$ equals
$$
\delta(n)=(\nu^o(n)-\nu^e(n))-(\lambda^o(n)-\lambda^e(n))+\Delta_3(n)-\Delta_3^{even}(n)
$$

and by (\ref{105}) and Lemma \ref{lem1}

\begin{equation}\label{107}
\delta(n)=\Delta_3^{odd}(n)+\varepsilon,
\end{equation}

where $|\varepsilon|\leq 2$.

Thus, by Theorem \ref{thm4} we have
\begin{equation}\label{108}
\lim_{n\rightarrow\infty}\frac{\ln \delta(n)}{\ln n}=\frac{\ln
3}{\ln4}.
\end{equation}

By the heuristic argument of the proportionality, we have
\begin{equation}\label{109}
\pi^o(n)-\pi^e(n)\approx \frac{3\pi(n)}{n}\delta(n).
\end{equation}

Now (\ref{108})-(\ref{109}) is equivalent to Conjecture 2.\;
$\blacksquare$

Table 2 compares on the powers of 4 the values of
$x(n)=\frac{\ln(\pi^o(n)-\pi^e(n))}{\ln n}$ and
$x^*(n)=\frac{\ln(\frac{3\pi(n)}{n}(\mu^e_3(n)-\mu^o_3(n))}{\ln n}$.

\begin{tab}

\begin{small}
$$
\begin{matrix}
&m &x(4^m) \;&x^*(4^m)  &m &x(4^m)\; &x^*(4^m)\\[8pt]
&2 &0.2500 \;&0.3962 \; &9\; &0.5983 \;&0.5974\\
&3 &0.3333 \;&0.4679 \; &10\; &0.6153 \;&0.6087\\
&4 &0.5574 \;&0.5109 \; &11\; &0.6237 \;&0.6186\\
&5 &0.5322 \;&0.5322 \; &12\; &0.6318 \;&0.6275\\
&6 &0.5736 \;&0.5537 \; &13\; &0.6364 \;&0.6354\\
&7 &0.5792 \;&0.5702 \; &14\; &0.6439 \;&0.6426
\end{matrix}
$$
\end{small}
\end{tab}

\section{On the increment of the excess of the odious
primes}

In conclusion let us consider the absolute value of the increment of
the excess of the numbers between the odious primes and the evil
primes on intervals $(0,2^n)$:
\begin{equation}\label{49}
\Delta(n)=\left|
(\pi^o(2^n)-\pi^e(2^n))-(\pi^o(2^{n-1})-\pi^e(2^{n-1}))\right|.
\end{equation}

By (\ref{107}),(\ref{109}), (\ref{49}) and Theorem \ref{1} we find
\begin{equation}\label{50}
\Delta(n)\approx
\begin{cases}3^{\frac{n-1}{2}}|\frac{\pi( 2^{n-1})}{2^{n-1}}-\frac{\pi(2^n)}{2^n}|,
\; if \; n\;is \;odd\\
3^{\frac n 2
-1}\left(3\frac{\pi(2^n)}{2^n}-\frac{\pi(2^{n-1})}{2^{n-1}}\right),
if \; n \;is \; even.
\end{cases}
\end{equation}

Notice that, by the Landau conjecture, $\pi(2n)\leq 2\pi(n), \;n\geq
3$ and therefore $\frac{\pi(2^{2n-1})}{2^{n-1}}\geq
\frac{\pi(2^n)}{2^n}, \; n\geq 2$. Unfortunately, this very
plausible conjecture was proved until now only for sufficiently
large $n$ \cite{RosScho}.

The following Table 4 illustrates the irregularity of the
distribution of $\Delta(n)$ (\ref{49}) in fact and by (\ref{50}) for
$n\geq 15$.

\begin{tab}
\begin{small}
$$
\begin{matrix}
&n &\Delta(n)\;&by(45)\\[8pt]
&15 &58 \;&19\\
&16 &492 \;&421\\
&17 &111 \;&42\\
&18 &1031 \;&1114\\
&19 &110 \;&98\\
&20 &3207 \;&2990\\
&21 &158 \;&238\\
&22 &8296 \;&8118\\
&23 &1416 \;&586\\
&24 &21790 \;&22229\\
&25 &1246 \;&1458\\
&26 &60294\;&61342\\
&27 &1570 \;&3707\\
&28 &170024 \;&170372
\end{matrix}
$$
\end{small}
\end{tab}

Notice that, although the phenomenon to a certain degree was
explained it remains very impressive that in spite of the ratio of
the numbers of primes in intervals $(2^{2t},
2^{2t+1}),\;(2^{2t-1},2^{2t})$is less than 2 but the value of
$\Delta$ (\ref{49}) for $t\geq 8$ more that $8, 9, 29, \ldots, 48,
108, \ldots$ times as large!

\itshape Conclusive remarks. \upshape\; 1)On the one hand, it is
interesting, using Theorem 1, to make the following steps  towards
justification of Conjecture 1. On the other hand, Conjecture 2 means
that the influence of the rest of the other steps in totality is
small. Nevertheless,  the full proof most likely requires more
strong  methods.

2)It is interesting to investigate the behavior of primes from the
considered point of view on the arithmetical progressions. For
example, on the progression $3t+2$ we expect on the whole an excess
of the \itshape evil \upshape primes since as one can show the
excess of the odd evil integers of the form $3t+2$ in interval
$[5,2^{2n-1})$ is equal to $3^{n-2}$, while on interval $[5,2^{2n})$
it is equal to 0. It is a topic for a separate article.

\end{document}